\documentclass[12pt]{article}
\usepackage{amssymb}
\newcommand{\bC}{{\mathbb C}}
\newcommand{\bZ}{{\mathbb Z}}
\newcommand{\bK}{{\mathbb K}}

\newcommand{\bF}{{\mathbf F}}
\newcommand{\bD}{{\mathbf D}}
\def\ord1{\mbox{\rm ord}}
\def\disc{\mbox{\rm disc}}
\def\mult{\mbox{\rm mult}}
\def\char{\mbox{\rm char}}
\newtheorem{thm}{Theorem}[section]
\newtheorem{lm}[thm]{Lemma}
\newtheorem{cor}[thm]{Corollary}

\newtheorem{rem}[thm]{Remark}

\title{A NOTE ON THE DISCRIMINANT
\footnotetext{
\begin{minipage}[t]{4.5in}
2010 Mathematics Subject Classification: 11C08, 13B25.\\
Key words and phrases: {\em discriminant of a polynomial,
nonsingular hypersurface}
\end{minipage}}}

\begin{document}
\author{by Arkadiusz P\l{}oski}
\maketitle

\begin{abstract}
Let $F(X,Y)=Y^d+a_1(X)Y^{d-1}+\dots+a_d(X)$ be a polynomial
in $n+1$ variables $(X,Y)=(X_1,\dots,X_n,Y)$ with coefficients
in an algebraically closed field $\bK$. Assuming that the
discriminant $D(X)=\disc_YF(X,Y)$ is nonzero we investigate
the order $\ord1_PD$ for $P\in\bK^n$. As application we get a discriminant
criterion for the hypersurface $F=0$ to be nonsingular.
\end{abstract}

\section*{Preliminaries}
Let $\bK$ be a fixed algebraically closed field of arbitrary characteristic.
A nonconstant polynomial $F$ in $N>0$ variables defines a hypersurface
$F=0$ which is by definition the set of all polynomials $aF$ where
$a\in\bK\setminus\{0\}$. Let $V(F)=\{P\in\bK^N:\,F(P)=0\}$
be the set of zeroes of the polynomial $F$. The order $\ord1_PF$ is the
lowest degree in the expansion of $F$ at $P$ obtained by translation
of coordinates. Thus $\ord1_PF>0$ if and only if $P\in V(F)$.
A point $P\in\bK^N$ is a simple (or nonsingular) point of the
hypersurface $F=0$ if $\ord1_PF=1$. If all points $P\in V(F)$ are
simple the hypersurface $F=0$ is called nonsingular.

For the notion of discriminant which is basic in this paper we refer the
reader to Abhyankar's book~\cite{A}, Lecture L4 (see also appendix IV
of~\cite{W} or~\cite{P}, \S~42).

\section{Result}

Let $F(X,Y)=Y^d+a_1(X)Y^{d-1}+\dots+a_d(X)\in\bK[X,Y]$ be a polynomial
in $n+1$ variables $(X,Y)=(X_1,\dots,X_n,Y)$ of degree $d>0$ in $Y$.
Let $D(X)=\mbox{disc}_YF(X,Y)$ be the $Y$-discriminant of $F$
(if $d=1$ then $D(X)=1$) and assume that $D(X)\neq 0$ in $\bK[X]$.
Then $F$ has no multiple factors in $\bK[X,Y]$.
Let $\pi:V(F)\rightarrow\bK^n$ be the projection given by
$\pi(a_1,\dots,a_n,b)=(a_1,\dots,a_n)$. In the sequel we write
$P=(a_1,\dots,a_n)$ and $Q=(a_1,\dots,a_n,b)$. We put
$\mult_Q\pi=\ord1_bF(P,Y)$ for $Q\in V(F)$. Then
$$
  \sum_{Q\in\pi^{-1}(P)}\mult_Q\pi=d
$$
for $P\in\bK^n$. Obviously we have $\#\pi^{-1}(P)\leq d$ with equality
if and only if $P\not\in V(D)$ (see~\cite{A}, Lecture~L4, Observation (0.2)).

The main result of this note is
\begin{thm}\label{thm11}
For any $P\in\bK^n$ we have $\ord1_PD\geq d-\#\pi^{-1}(P)$.
The equality $\ord1_PD=d-\#\pi^{-1}(P)$ holds if and only if
the following two conditions are satisfied
\begin{itemize}
\item[\rm(i)] $\ord1_QF=1$ for all $Q\in\pi^{-1}(P)$ i.e. all points
of $\pi^{-1}(P)$ are nonsingular,
\item[\rm(ii)] if $\char\bK=p$ then $p$ does not divide the multiplicities
$\mult_Q\pi$ for all $Q\in\pi^{-1}(P)$.
\end{itemize}
\end{thm}
We prove Theorem~\ref{thm11} in Section~3 of this note. Observe that
if $\char\bK=0$ then only Condition~(i) is relevant. Let us note a few
corollaries to Theorem~\ref{thm11}.
\begin{cor}\label{cor12}\mbox{}\\
If $\ord1_PD=0$ or $\ord1_PD=1$ then all points of $\pi^{-1}(P)$ are
nonsingular.
\end{cor}
Proof. If $\ord1_PD=0$ then $P\not\in V(D)$ and $\#\pi^{-1}(P)=d$ that is
$\ord1_PD=d-\#\pi^{-1}(P)=0$ and all points of $\pi^{-1}(P)$ are
nonsingular by Theorem~\ref{thm11}~(i). If $\ord1_PD=1$ then by the
first part of Theorem~\ref{thm11} we have $\#\pi^{-1}(P)\geq d-1$.
In fact $\#\pi^{-1}(P)=d-1$ because if $\#\pi^{-1}(P)=d$ then $P\not\in V(D)$
and $\ord1_PD=0$. We have then $\ord1_PD=d-\#\pi^{-1}(P)=1$ and again
by Theorem~\ref{thm11}~(i) we get the assertion. \rule{1ex}{1ex}
\begin{cor}\label{cor13} If $Q$ is a singular point of $F=0$
(i.e. $\ord1_QF>1$) then $\mult_Q\pi\leq\ord1_{\pi(Q)}D$.
\end{cor}
Proof. Let $P=\pi(Q)$. Since $Q\in\pi^{-1}(P)$ is singular we get by
Theorem~\ref{thm11} that $\ord1_PD>d-\#\pi^{-1}(P)$. On the other hand
it is easy to check that $\mult_Q\pi\leq d-\#\pi^{-1}(P)+1$. Therefore
we have $\mult_Q\pi\leq\ord1_{\pi(Q)}D$. \rule{1ex}{1ex}

\vspace{1ex}\noindent The following property of the discriminant is
well-known in the case $\bK=\bC$ (see \cite{W}, Appendix IV, Theorem~11~B).
\begin{cor}\label{cor14}
Let $\bF(A,Y)=Y^d+A_1Y^{d-1}+\dots+A_d\in\bZ[A,Y]$
be the polynomial with undeterminate coefficients $A=(A_1,\dots,A_d)$
and $\bD(A)=\disc_Y\bF(A,Y)$. Let $r(a)$ be the number of distinct roots
of the polynomial $\bF(a,Y)\in\bK[Y]$ where $a=(a_1,\dots,a_d)\in\bK^d$.
Then $\ord1_a\bD(A)\geq d-r(a)$.
The equality $\ord1_a\bD(A)=d-r(a)$ holds if $\char\bK=0$ or
$\char\bK=p$ and $p$ does not divide the multiplicities of roots of the
polynomial $\bF(a,Y)$.
\end{cor}
Proof. Observe that the hypersurface $\bF(A,Y)=0$ is nonsingular and
use Theorem~\ref{thm11}.
\begin{cor}\label{cor15}
{\rm(Discriminant Criterion for Nonsingular Hypersurfaces)}\\
Assume that $\bK$ is of characteristic zero. Then the hypersurface
defined by
$$
F(X,Y)=Y^d+a_1(X)Y^{d-1}+\dots+a_d(X)=0
$$
is nonsingular if and only if
$$
\ord1_P\bD(a_1(X),\dots,a_d(X))=\ord1_{a(P)}\bD(A_1,\dots,A_d)
$$
where $a(P)=(a_1(P),\dots,a_d(P))$, for all singular points $P$
of the discriminant hypersurface $D(X)=0$.
\end{cor}
Proof. Using Theorem~\ref{thm11} and Corollary~\ref{cor12} we see that
the hypersurface $F(X,Y)=0$ is nonsingular if and only if
$\ord1_PD(X)=d-\#\pi^{-1}(P)$ for all $P\in\bK^n$ such that $\ord1_PD(X)>1$.
On the other hand by Corollary~\ref{cor14} we have that
$d-\#\pi^{-1}(P)=\ord1_{a(P)}\bD(A_1,\dots,A_d)$. By definition of the
discriminant $D(X)=\bD(a_1(X),\dots,a_d(X))$ and the corollary follows.~\rule{1ex}{1ex}

\section{A property of the discriminant}

In the proof of Theorem~\ref{thm11} we need the following property of the
discriminant.
\begin{lm}\label{lm21}
Let $\bD(A_1,\dots,A_d)$ be the discriminant of the general polynomial
$\bF=Y^d+A_1Y^{d-1}+\dots+A_d$. Then
$$
  \bD(A_1,\dots,A_d)=(-1)^{d\choose 2}d^dA_d^{d-1}+\alpha_1(A_1,\dots,A_{d-1})A_d^{d-2}+\dots+\alpha_{d-1}(A_1,\dots,A_{d-1})
$$
in $\bZ[A_1,\dots,A_d]$, where $\ord1_0\alpha_i(A_1,\dots,A_{d-1})\geq i+1$
for $i=1,\dots,d-1$.
\end{lm}
Proof. Let
$$
I=\left\{(p_1,\dots,p_d)\in\bZ^d:\,p_i\geq 0\mbox{ for }i=1,\dots,d\mbox{ and }%
\sum_{i=1}^dip_i=d(d-1)\right\}\;.
$$
Then the discriminant $\bD(A_1,\dots,A_d)$ is equal to the sum of monomials of the form
$c_{p_1,\dots,p_d}A_1^{p_1}\dots A_d^{p_d}$ where $(p_1,\dots,p_d)\in I$
(see~\cite{A}, Lecture~L.4, Observation (03) or~\cite{P},\S~42).
It is easy to see that
\begin{itemize}
\item if $(p_1,\dots,p_d)\in I$ then $p_d\leq d-1$,
\item if $(p_1,\dots,p_d)\in I$ then $\sum_{i=1}^dp_i\geq d-1$
with equality if and only if $(p_1,\dots,p_d)=(0,\dots,0,d-1)$.
\end{itemize}
Therefore we get
$$
\bD(A_1,\dots,A_d)=c_{0,\dots,0,d}A_d^{d-1}+\alpha_1(A_1,\dots,A_{d-1})A_d^{d-2}+\dots+\alpha_{d-1}(A_1,\dots,A_{d-1})
$$
where $\ord1_0\alpha_i(A_1,\dots,A_{d-1})>i$ for $i=1,\dots,d-1$ and
$$
c_{0,\dots,0,d}=\bD(0,\dots,0,1)=\disc_Y(Y^d+1)=(-1)^{d\choose 2}d^d\quad\rule{1ex}{1ex}
$$
\begin{rem}\label{rem22}{\rm
Lemma~\ref{lm21} implies that $\ord1_0D(A_1,\dots,A_d)=d-1$.
Let $d(A_1,\dots,A_{d-1})=\disc_Y(Y^{d-1}+A_1Y^{d-2}+\dots+A_{d-1})$.
Then $\ord1_0d(A_1,\dots,A_{d-1})=d-2$. It is easy to check that
$\bD(A_1,\dots,A_{d-1},0)=\disc_Y(Y^d+A_1Y^{d-1}+\dots+A_{d-1}Y)=%
d(A_1,\dots,A_{d-1})A_{d-1}^2$. Therefore
$\ord1_0\alpha_{d-1}(A_1,\dots,A_{d-1})=\ord1_0\bD(A_1,\dots,A_{d-1},0)=%
\ord1_0d(A_1,\dots,A_{d-1})+2=(d-2)+2=d$.
} 
\end{rem}
From Lemma~\ref{lm21} and Remark~\ref{rem22} it follows that the Newton
diagram $\{(\ord1_0\alpha_k,d-k-1):\,\alpha_k\neq 0\}$ of the polynomial
$\bD=\alpha_0A_d^{d-1}+\alpha_1A_d^{d-2}+\dots+\alpha_{d-1}\in\bZ[A_1,\dots,A_{d-1}][A_d]$
intersects the axes at points $(0,d-1)$ and $(d,0)$. All remaining points
of the diagram lie strictly above the segment joining these two points.

\section{Proof}

Let $f(X,Y)\in\bK[[X,Y]]$ be a formal power series distinguished in $Y$ with
order $d>0$ i.e. such that $\ord1\,f(0,Y)=d$. Then by the Weierstrass
Preparation Theorem $f(X,Y)=F(X,Y)U(X,Y)$ in $\bK[[X,Y]]$ where
$F(X,Y)=Y^d+a_1(X)Y^{d-1}+\dots+a_d(X)\in\bK[[X]][Y]$ is a distinguished
polynomial and $U(0,0)\neq 0$.

Let $D(X)=\disc_YF(X,Y)$. Then we define $\tilde\mu(f)=\ord1_0D(X)-d+1$
if $D(X)\neq 0$ and $\tilde\mu(f)=+\infty$ if $D(X)=0$.
\begin{lm}\label{lm31}
We have $\tilde\mu(f)\geq 0$. The equality $\tilde\mu(f)=0$ holds
if and only if
\begin{itemize}
\item[\rm(i)] $\ord1_0f=1$,
\item[\rm(ii)] if $\char\bK=p$ then $p$ does not divide $d$.
\end{itemize}
\end{lm}
Proof. Let $D(X)\neq 0$. By Lemma~\ref{lm21} we have
$\ord1_0\bD(A_1,\dots,A_d)=d-1$ so 
$\ord1_0\bD(X)=\ord1_0\bD(a_1(X),\dots,a_d(X))\geq d-1$
since $\ord1_0a_i(X)\geq 1$ for $i=1,\dots,d$.
Again by Lemma~\ref{lm21} we can write
$$
  D(X)=(-1)^{d\choose 2}d^da_d(X)^{d-1}+\tilde D(X)
$$
where $\ord1_0\tilde D(X)>d-1$. Thus if $\ord1_0a_d(X)=1$ and $p$ does
not divide $d$ then $\ord1_0D(X)=d-1$. If $\ord1_0a_d(X)>1$ or $p$
divides $d$ then $\ord1_0D(X)>d-1$. This proves the lemma since
$\ord1_0F(X,Y)=1$ if and only if $\ord1_0a_d(X)=1$, and
$\ord1_0f=\ord1_0F$.~\rule{1ex}{1ex}
\begin{rem}\label{rem33}{\rm(see~\cite{T}, Lemma~5.10)
If $\char\bK=0$ then $\tilde\mu(f)$ is equal to the Milnor number of the
algebroid curve $f(cT,Y)=0$ where $X_i=c_iT$ ($i=1,\dots,r$) is a line
intersecting transversally the discriminant hypersurface $D(X)=0$.
} 
\end{rem}
\begin{lm}\label{lm32}
Let $F(X,Y)=Y^d+a_1(X)Y^{d-1}+\dots+a_d(X)\in\bK[X,Y]$ be a polynomial
such that $D(X)=\disc_YF(X,Y)\neq 0$. For any
$Q=(a_1,\dots,a_n,b)\in V(F)$ we put $F_Q(X,Y)=F(a_1+X,\dots,a_n+X_n,b+Y)$
(therefore $F_Q(X,Y)$ is distinguished in $Y$ with order $\mult_Q\pi$).
Then for every $P\in\bK^n$ we have
$$
\ord1_PD=\sum_{Q\in\pi^{-1}(P)}\tilde{\mu}(F_Q)+d-\#\pi^{-1}(P)\;.
$$
\end{lm}
Proof.Let $P=(a_1,\dots,a_n)$ and $r=\#\pi^{-1}(P)$.
Then $\pi^{-1}(P)=\{Q_1,\dots,Q_r\}$ where $Q_i=(P,b_i)$
with $b_i\neq b_j$ for $i\neq j$.
Let $d_i=\mult_{Q_i}\pi=\ord1_{b_i}F(a,Y)$.
Then $F(P,Y)=(Y-b_1)^{d_1}\dots(Y-b_r)^{d_r}$.
We have to check that
\begin{equation}\label{eq1}
\ord1_PD=\sum_{i=1}^r\tilde\mu(F_{(P,b_i)})+d-r\;.
\end{equation}
Let $F_P(X,Y)=F(a_1+X_1,\dots,a_n+X_n,Y)$ and
$D_P(X)=D(a_1+X_1,\dots,a_n+X_n)$. Then $D_P(X)=\disc_YF_P(X,Y)$
and $\ord1_PD(X)=\ord1_0D(X)$. Moreover
$(F_P)_{(0,b_i)}=F_{(P,b_i)}$ for $i=1,\dots,r$ and it suffices to prove
(\ref{eq1}) for $F_P$ at $0\in\bK^n$. Henceforth we assume that
$P=0\in\bK^n$. First let us consider the case $r=1$. Then (\ref{eq1})
reduces (for $P=0$) to the formula
\begin{equation}\label{eq2}
\ord1_0D=\tilde\mu(F_{(0,b_1)})+d-1\mbox{ provided that }F(0,Y)=(Y-b_1)^d\;.
\end{equation}
The polynomials
$F_{(0,b_1)}(X,Y)=F(X,b_1+Y)$ and $F(X,Y)$ have the same $Y$-discriminant,
hence (\ref{eq2}) follows directly from the definition of $\tilde\mu$.
Suppose that $r>1$. By Hensel's Lemma we have
$$
  F(X,Y)=\prod_{i=1}^rF_i(X,Y)\mbox{ in }\bK[[X]][Y]\mbox{ with }
  F(0,Y)=(Y-b_i)^{d_i}\;.
$$
Let $D_i(X)=\disc_YF_i(X,Y)$ for $i=1,\dots,r$ and
$R_{ij}(X)=Y$-resultant $F_i(X,Y),F_j(X,Y)$ for $i\neq j$.
By the product formula for the discriminant
$$
  D(X)=\prod_{i=1}^rD_i(X)\prod_{1\leq i<j\leq r}R_{ij}(X)^2
$$
we get
$$
  \ord1_0D(X)=\sum_{i=1}^r\ord1_0D_i(X)
$$
since $R_{ij}(0)=(b_i-b_j)^{d_id_j}\neq 0$ for $i\neq j$.

By Formula (\ref{eq2}) (which applies to the polynomials with coefficients
in $\bK[[X]]$) we have
$$
  \ord1_0D_i=\tilde\mu((F_i)_{(0,b_i)})+d_i-1\mbox{ for }i=1,\dots,r\;.
$$
Since $(F_i)_{(0,b_i)}$ and $F_{(0,b_i)}$ are associated in $\bK[[X,Y]]$
we can write
$$
  \ord1_0D_i=\tilde\mu(F_{(0,b_i)})+d_i-1
$$
and
$$
  \ord1_0D=\sum_{i=1}^r(\tilde\mu(F_{(0,b_i)})+d_i-1)=%
  \sum_{i=1}^r\tilde\mu(F_{(0,b_i)})+d-r\;.~\rule{1ex}{1ex}
$$

\vspace{1ex}\noindent{\bf Proof of Theorem~\ref{thm11}}\\
Theorem~\ref{thm11} follows directly from Lemmas~\ref{lm31} and~\ref{lm32}.~\rule{1ex}{1ex}

\medskip
\noindent {\small
Chair of Mathematics\\
Technical University\\
Al. 1000 L PP7\\
25-314 Kielce, Poland\\
e-mail address: matap@tu.kielce.pl}


\begin{thebibliography}{G}

\bibitem[1] {A}
S.S. Abhyankar, {\em Lectures on Algebra I\/},
World Scientific, 2006.

\bibitem[2] {P}
O.~Perron,
{\em Algebra I\/},
Springer, Berlin, 1927.

\bibitem[3] {T}
B.~Teissier, {\em Cycles \'{e}vanescences, sectiones planes et
conditions de Whitney\/},
Ast\'{e}risque 1973, 7--8, 285--362.

\bibitem[4] {W}
H.~Whitney, {\em Complex Analytic Varieties\/},
Addison--Wesley 1972.

\end{thebibliography}
\end{document}